\documentstyle[12pt,twoside]{amsart}

\renewcommand{\c}[0]{{\mathbb C}}  

\renewcommand{\o}[0]{{\cal O}} 
\newcommand{\z}[0]{{\mathbb Z}}

\renewcommand{\r}[0]{{\mathbb R}}

\newcommand{\p}[0]{{\mathbb P}}

\newcommand{\q}[0]{{\mathbb Q}}
\newcommand{\map}[0]{\dasharrow}
\newcommand{\qtq}[1]{\quad\mbox{#1}\quad}
\newcommand{\spec}[0]{\operatorname{Spec}}
\newcommand{\pic}[0]{\operatorname{Pic}}

\newcommand{\proj}[0]{\operatorname{Proj}}

\newcommand{\nec}[1]{\overline{NE}({#1})}

\newsymbol\subsetneq 2328
\newsymbol\onto 1310
\def\into{\DOTSB\lhook\joinrel\rightarrow}
\newsymbol\twoheadrightarrow 1310

\newtheorem{thm}{Theorem}[section]

\newtheorem{lem}[thm]{Lemma}
\newtheorem{cor}[thm]{Corollary}

\newtheorem{prop}[thm]{Proposition}

\theoremstyle{definition}
\newtheorem{defn}[thm]{Definition}

\newtheorem{say}[thm]{}

\newtheorem{exrc}[thm]{Exercise}

\newtheorem{ack}{Acknowledgments}        
\newtheorem{notation}[thm]{Notation}

\theoremstyle{remark}

\begin{document}
\bibliographystyle{amsplain}

\title{Real Algebraic Surfaces}
\author{J\'anos Koll\'ar}

\maketitle

The aim of these lectures is to study real algebraic surfaces  with the
help of the minimal model program. 
I  mainly concentrate on those surfaces $F$ which are rational over
$\c$. In some sense this is the class where the theory over $\c$ and
over $\r$ differ the most.

The study of real algebraic curves and surfaces is very old, in fact it
predates algebraic geometry over $\c$. From our point of view the first
significant result was the classification of real cubic surfaces by
\cite{Schlafli1863}. \cite{Zeuthen1874} studied plane curves of degree
4 and their bitangents, which is essentially equivalent to degree 2 Del
Pezzo surfaces (though he related them to cubic surfaces only). 
The first systematic study was done in \cite{Comessatti13} and
\cite{Comessatti14}. These papers presented a   breakthrough whose
extent has not been appreciated for over 60 years.
(For instance, the main theorem of \cite{Comessatti14}  disproves a
conjecture of \cite{Nash52}.) The monograph
\cite{Segre42} treats cubic surfaces in detail. 

The development of the theory of surfaces over nonclosed fields was
aimed primarily at arithmetic questions, but the results of
\cite{Segre51, Manin66, Iskovskikh80} are very interesting even in the
real case.

A systematic study of real algebraic surfaces is given in
\cite{Silhol89}. This book develops the theory from the point of view
of topology, Galois cohomology and minimal model theory.

The aim of my lectures is to  investigate 
real algebraic surfaces  using primarily the theory of minimal models.
For surfaces this is    mostly a question of
personal preference, however, this approach becomes significant for 
real algebraic  threefolds.

The  topological methods of real algebraic geometry are basically
homological. That is, they give information about the homology
of the set of real points. In the surface case the set of real points
is a topological manifold of real dimension 2. Compact manifolds 
of real dimension 2 are classsified by their homology, hence in this
case there is no need for more detailed information.
By contrast, for compact manifolds 
of real dimension 3  their homology carries relatively little
information. The geometric methods of these lectures can be generalized
to threefolds. Thus these notes form an introduction to my papers
\cite{rat}, though the 3--dimensional methods and results are never
mentioned explicitly.

\begin{notation} If $X$ is a scheme over $\r$ then $X(\r)$ denotes its
set of real points with the Euclidean topology and
$X_{\c}:=X\times_{\spec \r}\spec \c$ denotes the complexification of 
$X$.  Sometimes  I write $X_{\r}$ instead of $X$ to emphasize that we
are over $\r$. One should keep in mind the distinction between
$X_{\r}$ and $X(\r)$. The first is a scheme over $\r$, thus among its
closed points we find all conjugate pairs of complex points of
$X_{\c}$, but $X(\r)$ is the collection of the real points only.

$\p^n$ stands for projective $n$-space  as a  scheme over $\q$.
Thus $\p^n_{\r}$ is  projective $n$-space  as a  scheme over $\r$
and $\p^n(\r)$ is the set of real points, customarily denoted by
$\r\p^n$. 

When talking about a point, curve or vector bundle etc. on a real variety
$X$, these are, by definiton, objects defined over $\r$. I will always
explicitly mention if we are considering 
a point, curve or vector bundle etc. on $X_{\c}$. 

If $Z\subset X_{\c}$ is a subscheme, then $\bar Z\subset X_{\c}$
denotes the subscheme defined by conjugate equations. Then $Z+\bar
Z\subset X_{\c}$ can be  defined by real equations and it can be viewed
as a subscheme of $X_{\r}$. That is, there is a unique subscheme
$Z_{\r}\subset X_{\r}$ such that $Z_{\r}\times_{\spec \r}\spec
\c=Z+\bar Z$. 

A property of $X$ (irreducible, normal, smooth etc.) always refers to
the scheme theoretic property of $X$. Thus if $X$ is smooth then it is
smooth not just at all points of $X(\r)$, but also at all complex points.
When talking about irreducibility (normality etc.)  of $X_{\c}$, I say
$X$ is  ``geometrically irreducible" or ``irreducible over $\c$". In some
cases this does not matter (for instance $X$ is normal or smooth iff
$X_{\c}$ is) but in other cases the two versions are different ($X$ may
be irreducible but not geometrically irreducible).

The logical continuation of this terminology is to say that two
irreducible 
$\r$-schemes $X,Y$ are birational iff there is a map $f:X\map Y$ defined
over $\r$ which is birational. Similarly, $X$ is rational iff it is
birational to $\p^n_{\r}$ over $\r$. Unfortunately many authors use the
word ``rational" to mean ``geometrically rational". To avoid confusion,
I  will use the expressions ``rational over $\r$" and ``rational over
$\c$".
\end{notation}

\section{Minimal models of real algebraic surfaces}

\begin{defn} Let $X$ be a   variety over a field $k$.
A {\it 1--cycle} on $X$ is a formal linear combination $C=\sum c_iC_i$,
where the $C_i\subset X$ are irreducible, reduced  and proper curves. A
1--cycle is called {\it effective} if $c_i\geq 0$ for every $i$.

Two 1--cycles $C,C'$ are {\it numerically equivalent} if
$(C\cdot D)=(C'\cdot D)$ for every Cartier divisor $D$ on $X$.
1--cycles with real coefficients modulo numerical equivalence form a
vectorspace, denoted by $N_1(X)$. $N_1(X)$ is finite dimensional by the
Theorem of the base of N\'eron--Severi (cf.
\cite[p.447]{Hartshorne77}).  Its dimension, denoted by $\rho(X)$, is
called the {\it Picard number} of $X$.  

Effective 1--cycles generate a cone  $NE(X)\subset N_1(X)$. Its closure
in the Euclidean topology $\overline{NE}(X)\subset N_1(X)$ is called
the {\it cone of curves} of $X$.

If $K_X$ is Cartier (or at least some multiple of $K_X$ is Cartier) then 
set 
$$
\overline{NE}(X)_{K\geq 0}:=\{z\in \overline{NE}(X)\vert (z\cdot
K_X)\geq 0\}.
$$

Let $V\subset \r^n$ be a closed convex cone. For $v\in V$, a ray
$\r^{\geq 0}v\subset V$ is called {\it extremal} if $u,u'\in V$, $u+u'\in
\r^{\geq 0}v$ implies that $u,u'\in\r^{\geq 0}v$. Intuitively: $\r^{\geq 0}v$
is an edge of
$V$.

An extremal ray $\r^{\geq 0}z\subset \overline{NE}(X)$ is called {\it
$K_X$-negative} if $(z\cdot K_X)<0$. This does not depend on the choice
of $z$ in the ray.

Let $R\subset \overline{NE}(X)$ be a  ray. A {\it
contraction} of
$R$ is a morphism $f_R:X\to X'$ such that $(f_R)_*\o_X=\o_{X'}$ and a
curve $C\subset X$ is mapped to a point iff $[C]\in R$.
\end{defn}

\begin{exrc}\label{er.noncontr.ex}
  Show that the contraction of a  ray is unique
(if it exists).  Also, if $X'$ is projective then $R$ is an extremal
ray. Find examples of extremal rays which can not be contracted.
\end{exrc}

\begin{exrc}\label{int.cone.exrc}
  Let $F$ be a smooth projective surface and $z\in 
\overline{NE}(F)$ a 1--cycle such that $(z^2)>0$. 
Then $z$ or $-z$ is in the interior of $\overline{NE}(F)$.
Thus if $\r^{\geq 0}z$ is extremal and $(z^2)>0$ then 
$\overline{NE}(F)$ is 1--dimensional.
\end{exrc}

We use the following description of the cone of curves of  smooth
surfaces over $\c$. The result is essentially equivalent to the theory
of minimal models of surfaces developed around the turn of the 20th
century. This formulation (and its higher dimensional generalization)
is due to \cite{Mori82}. See also \cite{koll96, km98} for proofs.

\begin{thm}[Cone Theorem]\label{cone.thm.c}
 Let $F$ be a smooth projective surface over an
algebraically closed field. Then there are curves $C_i\subset F$
such that
$$
\overline{NE}(F)=\overline{NE}(F)_{K\geq 0}+\sum \r^{\geq 0}[C_i],
$$
and  the $\r^{\geq 0}[C_i]$ are $K_F$-negative extremal rays of 
$\overline{NE}(F)$. 

Moreover, we can assume that each 
$C_i\subset F$ is a smooth rational curve and
$(C_i^2)\in\{-1,0,1\}$. If $(C_i^2)=1$ (resp. $(C_i^2)=0$) for
some $i$  then $F\cong
\p^2$ (resp. $F$ is a minimal ruled surface over a curve). \qed
\end{thm}

Let now $F$ be a smooth projective surface over $\r$. If $C$ is a
1--cycle on  $F_{\c}$   then $C+\bar C$ is a 1--cycle  on $F$, and
every 1--cycle  on $F$ arises this way, at least if we use rational or
real coefficients. Thus (\ref{cone.thm.c}) immediately gives:

\begin{thm}[Cone Theorem over $\r$]\label{cone.thm.r}
 Let $F$ be a smooth projective surface over $\r$. Then there are 
smooth rational curves
$C_i\subset F_{\c}$ with $(C_i^2)\in\{-1,0,1\}$ such that 
$$
\overline{NE}(F)=\overline{NE}(F)_{K\geq 0}+\sum \r^{\geq 0}[C_i+\bar C_i],
$$
and  the $\r^{\geq 0}[C_i+\bar C_i]$ are $K_F$-negative extremal rays of 
$\overline{NE}(F)$. 
\end{thm}

Proof. There is  one point that we need to be careful about.
Namely, it happens frequently that $C_i$ gives an extremal ray but
$C_i+\bar C_i$ does not. So we have to throw  away some of the $C_i$
appearing in (\ref{cone.thm.c}). \qed
\medskip

\begin{say}[Geometric irreducibility]  Let $X\subset \p^n$ be a variety
over
$\c$ and $\bar X$ the variety defined by conjugate equations. The
disjoint union of $X$ and $\bar X$ is invariant under conjugation, and
so there is a real variety $Y_{\r}$ such that $Y_{\c}\cong X\cup \bar
X$. Such real varieties are not particularly interesting since the
theory of $Y_{\r}$ over $\r$ is equivalent to the theory
of $X$ over $\c$. Thus it is reasonable to restrict our attention to
real varieties $Y$ such that $Y_{\c}$ is irreducible, that is, $Y$ is
geometrically irreducible.

Of course, during a proof we may run into a subvariety of $Y_{\r}$
which is geometrically reducible, and these have to be dealt with
appropriately. Thus we can not ignore such varieties completely.
\end{say}

\begin{defn} Let $S$ be a  smooth projective surface  over  a field $k$.
$S$ is called a {\it Del Pezzo surface} if $S$ is geometrically
irreducible and
$-K_S$ is ample. It is called {\it minimal} (over $k$) if $\rho(S)=1$.

$S$, together with a morphism to a smooth curve $f:S\to B$ is called a
{\it conic bundle} if every fiber is isomorphic to a plane conic.
A conic bundle is called {\it minimal} if $\rho(S)=2$.
\end{defn}

The geometric description and meaning of the extremal rays occurring in
(\ref{cone.thm.r}) is given in the next result:

\begin{thm}\label{mmp.over.R}
 Let $F$ be a smooth projective geometrically irreducible surface  over
$\r$ and $R\subset \nec{F}$ a $K_F$-negative extremal ray. Then $R$ can
be contracted
$f:F\to F'$, and we obtain  one of the following cases:
\begin{enumerate}
\item[(B)] (Birational)   $F'$ is a smooth projective surface  over
$\r$  and $\rho(F')=\rho(F)-1$.
 $F$ is the blow up of $F'$ at  a closed point $P$. We have two
cases:
\begin{enumerate}
\item  $P\in F'(\r)$, or
\item  $P$ is a pair of conjugate points.
\end{enumerate}

\item[(C)] (Conic bundle) $B:=F'$ is a smooth curve, $\rho(F)=2$ and
$F\to B$ is a conic bundle. The fibers   $f^{-1}(P): P\in B(\c)$ are
smooth, except for an even number of  ponts $P_1,\dots,P_{2m}\in B(\r)$. 
$(K_F^2)=8(1-g(B))-2m$.

\item[(D)] (Del Pezzo surface) $F'$ is a point, $\rho(F)=1$, $-K_F$ is
ample and we have one of the following cases
\begin{enumerate}
\item  $(K_F^2)=9$ and $F\cong  \p^2$.
\item  $(K_F^2)=8$ and  $F\cong (x_0^2+x_1^2+x_2^2-x_3^2=0)\subset
\p^3$.
\item   $(K_F^2)=2$.
\item   $(K_F^2)=1$.
\end{enumerate}
\end{enumerate}
\end{thm}

Proof. By (\ref{mmp.over.R}), 
there is a   curve $C\cong \p^1$  over $\c$  such that 
  $C+\bar
C$ generates $R$  and $(C^2)\in \{-1,0,1\}$.

 We consider various possibilities.

Assume first that $(C^2)=-1$. If $C=\bar C$ then 
the contraction of $C$ in $F_{\c}$ is defined over $\r$, thus
$F$
is the blow up of a surface at a real point.
(This is Castelnuovo's contraction theorem, cf.
\cite[V.5.7]{Hartshorne77}.) If $C$ and $\bar C$ are disjoint, then we
can contract them simultaneously over
$\r$ to obtain
$f:F\to F'$ which is an isomorphism near $F(\r)$.

If $(C^2)=-1$ and $(C\cdot \bar C)=1$ or $(C^2)= (C\cdot \bar
C)=0$, then $C+\bar C$ has selfintersection 0.
From Riemann--Roch we obtain that 
$$
h^0(F, \o_F(m(C+\bar C))
\geq \chi (F, \o_F(m(C+\bar C))= m+\chi(\o_F),
$$
thus $m(C+\bar C)$ moves in a linear system for $m\gg 1$. It's moving
part is base point free by (\ref{bpf.exer}).
Let $f:F\to B$ be the Stein factorization of the resulting morphism.

Let  $A\subset F$ be an irreducible    fiber. If $A$ is a multiple
fiber, write it as $A=mA_1$.  Since $[A_1]\in R$,
$$
2g(A_1)-2=(A_1^2)+(A_1\cdot K_F)=(A_1\cdot K_F)<0.
$$
 Thus $A_1$ is  isomorphic to a smooth conic over $\r$ (\ref{conics.r})
and $(A_1\cdot K_F)=-2$. The generic fiber $A_g$ is not multiple, so
$(mA_1\cdot K_F)=(A_g\cdot K_F)=-2$ which shows that there are no
multiple fibers.

Let $A_1+ A_2=f^{-1}(b)$ be a  reducible   fiber over $\c$, where $A_1$
is an irreducible and reduced curve.  In particular,  $(A_1^2)<0$. 
$\bar A_1+ \bar A_2=f^{-1}(\bar b)$ is also a fiber.
If $b\neq \bar b$ then  $A_1$ is disjoint from $\bar A_1$, thus
$((A_1+\bar A_1)^2)=(A_1^2)+(\bar A_1^2)<0$.
$[A_1+\bar A_1]+[A_2+\bar A_2]\in R$ and 
$[A_1+\bar A_1]\not\in R$, a contradiction. Thus all singular fibers
lie over real points $P_1,\dots, P_r\in B(\r)$.

Therefore, $A_1+\bar A_1\subset f^{-1}(b)$. Every fiber of $f$ is
irreducible over $\r$, thus $A_1+\bar A_1= f^{-1}(b)$.
We get that $(A_1\cdot K_X)=-1$ and so $A_1$ is a $-1$-curve
and $A_1+\bar A_1$  is isomorphic to a pair of conjugate lines in
$\p^2$. 
$F$ is a conic bundle over $B$ by (\ref{conic.bund}).

Over $\c$ we can contract one of the components of every singular fiber
to obtain a minimal ruled surface.  The selfintersection  number of the
canonical class of a minimal ruled surface is $8(1-g(B))$, and each
singular fiber drops this number by 1. We see in (\ref{mmp.over.R.top})
that the number of singular fibers is even.

Assume next that $(C^2)=-1$ and
 $r:=(C\cdot \bar C)\geq 2$. Then 
$$
((C+\bar C)^2)=-2+2(C\cdot \bar C)=2r-2>0.
$$
By (\ref{int.cone.exrc}) 
this implies that $\nec{F}$  is 1-dimensional, hence $-K_F\equiv
a(C+\bar C)$ for some $a>0$.
$(-K_F\cdot (C+ \bar C))=2$, thus 
$$
(C+\bar C)\equiv (1-r)K_F\qtq{and} 2r-2=(1-r)^2(K_F^2).
$$
This gives the possibilities  $r=2, (K_F^2)=2$ or $r=3, (K_F^2)=1$.

If $(C^2)=0$ and
 $r:=(C\cdot \bar C)\geq 1$ then a computation as above gives that
$8=r(K_F^2)$, which allows too many cases. It is better to consider
this geometrically.

By (\ref{cone.thm.c}), 
$C$ is a fiber of  a $\p^1$-bundle $g:F_{\c}\to D$ over $\c$
and $\bar C$ is a (possibly multiple) section of $g$.
Thus $D$ is rational. By the classification of minimal ruled surfaces,
either $F_{\c}\cong \p^1\times \p^1$ and we are done by
(\ref{quadric.lem}), or
$g$ has a unique section
$E$ with negative selfintersection. $E$ is then defined over $\r$, thus
$\rho(F)=2$, a contradiction.

We are left with the case when   $F_{\c}\cong \c\p^2$
and $C$ is a line in $\c\p^2$. Then $\bar C$ is another line
and $C$ and $\bar C$ intersect in a unique point, which is therefore
real.  We can get another real point, and so also a real line. Thus
 $\o_F(1)$ is defined over
$\r$ and 
$F\cong
\r\p^2$. \qed
\medskip

As  a consequence we obtain the minimal model program (MMP for short)
for real algebraic surfaces:

\begin{thm}[MMP for surfaces]\label{mmp.surf.thm}
Let $F$ be a smooth projective geometrically irreducible surface over
$\r$. Then there is a sequence of morphisms
$$
F=F_0\stackrel{f_0}{\to} F_1 \to \cdots F_{m-1}\stackrel{f_{m-1}}{\to}
F_m=F^*
$$
such that each $f_i:F_i{\to} F_{i+1}$ is a birational contraction as in
(\ref{mmp.over.R}.B) and $F^*$ satisfies precisely one of the following
properties:
\begin{enumerate}
\item[(M)] (Minimal model)  $K_{F^*}$ is nef. (That is, it has
nonnegative intersection number with every curve in $F^*$.)

\item[(C)] (Conic bundle)   $F^*$ is a conic bundle
over a curve 
$f:F^*\to B$. In particular,  $\rho(F^*)=2$

\item[(D)] (Del Pezzo surface)  $\rho(F^*)=1$, $-K_{F^*}$ is ample
and $F^*$ is among those listed in (\ref{mmp.over.R}.D).
\end{enumerate}
\end{thm}

Proof. We do the steps of (\ref{mmp.over.R}.B) as long as we can.
$\rho(F_{i+1})=\rho(F_i)-1$, so eventually we reach $F^*=F_m$ where we
can not perform a contraction as in (\ref{mmp.over.R}.B). 
If $K_{F^*}$ is nef then we have a minimal model.

If $K_{F^*}$ is not nef, then by (\ref{mmp.over.R}) we can perform  a
contraction as in (\ref{mmp.over.R}.C--D). This gives our last two
cases.\qed

\begin{exrc}\label{conics.r} Let $k$ be a field and $C$ a smooth
 projective curve over $k$ such that $C_{\bar k}\cong \p^1$. Show that
$C$ is  isomorphic to a smooth conic over $k$. Also, $C\cong \p^1_k$
iff $C(k)\neq\emptyset$.
\end{exrc}

\begin{exrc}\label{conic.bund} Let $S$ be a smooth projective surface
over $\c$ and $f:S\to B$ a morphism to a smooth curve. Assume that
$f^{-1}(b)\cong \p^1$ for  some $b\in B$ and every fiber has at most 2
irreducible components. Show that 
$-K_F$ is very ample on the fibers, $f_*\o_F(-K_F)$ is  a rank 3
vector bundle over $B$ and we have an injection
$$
F\into \proj_B f_*\o_F(-K_F).
$$
Under this injection the fibers of $f$ become conics. Such a surface is
called
  a {\it conic bundle} over $B$.
\end{exrc}

\begin{exrc}\label{projspace.ex}
 Let $X$ be a variety over $\r$ such that $X_{\c}\cong
\p^n$.  Show that $X\cong \p^n$ if $n$ is even, but not necessarily if
$n$ is odd.
\end{exrc}

\begin{exrc}[Cohomology commutes with base change]\label{coh.b.c}
 Let $X$ be a variety over $\r$  and $F$ a coherent sheaf on $X$. Show
that
$$
H^i(X,F)\otimes_{\r}\c\cong H^i(X_{\c}, F_{\c}).
$$
\end{exrc}

\begin{exrc}\label{bpf.exer} Let $F$ be a smooth projective surface and
$D\subset F$ an irreducible curve such that $(D^2)=0$. Then the moving
part of $|mD|$ is either empty or base point free.
\end{exrc}

\begin{notation} $Q^{r,s}$ denotes the quadric hypersurface
$(x_1^2+\cdots +x_r^2-x_{r+1}^2-\cdots-x_{r+s}^2=0)$.
$Q^{2,1}$ is isomorphic to $\p^1$ by (\ref{conics.r}).
\end{notation}

\begin{lem}\label{quadric.lem}
 Let $F$ be a smooth projective surface over $\r$ such that
$F_{\c}\cong \p^1\times \p^1$. Then
 one of the
following holds:
\begin{enumerate}
\item $F\cong Q^{2,2}\cong Q^{2,1}\times Q^{2,1}$, $\rho(F)=2$ and
$F(\r)\sim S^1\times S^1$,
\item $F\cong Q^{3,1}$, $\rho(F)=1$ and $F(\r)\sim S^2$,  
\item $F\cong Q^{4,0}\cong Q^{3,0}\times Q^{3,0}$, $\rho(F)=2$ and
$F(\r)=\emptyset$,
\item $F\cong Q^{3,0}\times \p^1$, $\rho(F)=2$  and
$F(\r)= \emptyset$.
\end{enumerate}
\end{lem}

Proof. Let $C\subset F_{\c}$ be one of the rulings. Then $\bar C$ is
another ruling, thus either $(C\cdot \bar C)=0$ or $(C\cdot \bar C)=1$.

If $(C\cdot \bar C)=0$ then the linear system $|C+\bar C|$ is defined
over $\r$ and maps $F$ onto a conic. Similarly for the other rulings,
thus $F$ is the product of two conics.
All 3 possibilities are listed.

If $(C\cdot \bar C)=1$ then $\o_F(C+\bar C)$ is a line bundle on $F$
which is of type 
$\o_{F_{\c}}(1,1)$ over $\c$. Thus its global sections embed $F$ as a
quadric.
$Q^{3,1}$ is the only quadric not yet accounted for.\qed
\medskip

In the above proofs we had to establish several times that certain line
bundles on $F_{\c}$ are defined over $\r$. This is frequently a quite
subtle point. Some aspects of it are treated in the next exercise.

\begin{exrc}\label{pic=gal.inv}  Let $X$ be a scheme over $\r$ and
$L$  a line bundle on $X$. Then $L_{\c}$ is a line bundle on $X_{\c}$
and $L_{\c}\cong \bar L_{\c}$.
Thus if  $M$ is a line bundle on $X_{\c}$ and $M\not\cong \bar M$,
then $M$ is not the complexification of a real line bundle.

Find a curve $C$ over $\r$ and a line bundle $M$ on $C_{\c}$ such that
$M\cong \bar M$ but $M$ is not the complexification of a real line
bundle.

Let $X$ be a scheme over $\r$ and  $M$ a line bundle on $X_{\c}$ such
that
$M\cong \bar M$. Show that $M^{\otimes 2}$ is the  complexification of a
real line bundle. If $X$ is connected, reduced and $X(\r)\neq \emptyset$
then $M$ itself is the complexification of a real line
bundle.

 More generally, let   $X_K$  be an integral scheme defined
over a field
$K$ and 
$L\supset K$ a Galois extension with Galois group $G$. 
Show that if $X_K$ 
has a $K$-point then $\pic (X_K)=\pic (X_L)^G$.
\end{exrc}

\section{The topology of  $F(\r)$}

In this section we study the MMP from the topological point of view.
The main results of this section are already in \cite{Comessatti14}.

\begin{notation} $M\uplus N$ denotes the disjoint union of $M$ and $N$.
$\uplus rN$ denotes the disjoint union of $r$ copies of  $N$.
$M\# N$ denotes the connected sum of two manifolds $M$ and $N$ (which
are assumed to have the same dimension).
$\# rN$ denotes the connected sum  of $r$ copies of  $N$. 
(By definition, $\#0M=S^{\dim M}$.) $M\sim N$ denotes
that $M$ and $N$ are homeomorphic.
\end{notation}

One can give a complete topological description of the 
various contractions in (\ref{mmp.over.R}):

\begin{thm}\label{mmp.over.R.top}
 Let $F$ be a smooth projective geometrically irreducible surface  over
$\r$ and $R\subset \nec{F}$ a $K_F$-negative extremal ray. The  following
is the topological description of the corresponding contraction:
\begin{enumerate}
\item[(B)] (Birational)   $F$ is the blow up of $F'$ at  a closed point
$P$. We have two cases:
\begin{enumerate}
\item  If $P\in F'(\r)$ then  $F(\r)\sim F'(\r)\#\r\p^2$.
\item  If $P$ is a pair of conjugate points then $F(\r)\sim F'(\r)$.
\end{enumerate}

\item[(C)] (Conic bundle) $f:F\to B $  is a conic bundle with singular
fibers  
$f^{-1}(P_1),\dots,f^{-1}(P_{2m})$. Then
$$
F(\r)\sim \uplus mS^2\uplus N_1\uplus\cdots \uplus N_b,
$$
where $b$ is the number of connected components of $B(\r)$ 
which do not contain any of the points  $P_i$ and each $N_i$ is either a
torus or a Klein bottle.

\item[(D)] (Del Pezzo surface) There are 4 cases:
\begin{enumerate}
\item  If $(K_F^2)=9$ then $F(\r)\sim  \r\p^2$.
\item  If $(K_F^2)=8$ then $F(\r)\sim S^2$.
\item  If  $(K_F^2)=2$ then $F(\r)\sim \uplus 4S^2$.
\item   If $(K_F^2)=1$ then $F(\r)\sim \r\p^2\uplus 4S^2$.
\end{enumerate}
\end{enumerate}
\end{thm}

Proof.  Blowing up replaces a point with all tangent directions 
through that point. So we remove a disc and put in an interval bundle
over $S^1$ whose boundary is connected. This is a M\"obius strip and so
 $F(\r)\sim
F'(\r)\#\r\p^2$. 

In the conic bundle case, 
let $M\sim S^1$ be a connected component of $B(\r)$. If none of the
$P_i$ lie on $M$ then $F(\r)\to B(\r)$ is  a smooth $S^1$-bundle over
$M$, this gives either a
torus or a Klein bottle. If $k$ of the points 
  $P_1,\dots,P_k\in M\sim S^1$ correspond to singular
fibers then, after reindexing, they divide $M$ into $k$ intervals
$[P_i,P_{i+1}]$ (subscript
$\mod k$).   $F(\r)$ is alternatingly empty or a copy of $S^2$ over
the intervals. Thus $k$ is even.

In the Del Pezzo case we are done if $(K_F^2)=9,8$. The cases
$(K_F^2)=2,1$ are considerably harder. They follow from
(\ref{deg2.equiv.thm}) and (\ref{deg1.equiv.thm}).\qed
 \medskip

Using (\ref{mmp.over.R.top}) it is easy  to determine which 
2-manifolds occur as $F(\r)$ for geometrically rational surfaces $F$.
The   conclusion is that orientable surfaces of genus $>1$ do not
occur. This is the main result of 
\cite{Comessatti14}.

\begin{thm}\label{c-rat.top}
 Let $F$ be a smooth,
projective surface  over
$\r$ such that $F_{\c}$ is rational.
Then  one of the following holds:
\begin{enumerate}\setcounter{enumi}{-1}
\item    $F(\r)=\emptyset$.
\item   $F(\r)\sim S^1\times S^1$.
\item   $F(\r)\sim \#r_1\r\p^2\uplus \cdots \uplus \#r_m\r\p^2$ for some
$r_1,\dots,r_m\geq 0$.  
\end{enumerate}
All these cases do occur.
\end{thm}

Proof. Apply the MMP over $\r$ to get $F=F_1\to F_2\to\cdots$.
We prove the theorem by induction on the number of blow ups in the
sequence.  If $F_i\to F_{i+1}$ is the inverse of the blowing up
of a real point, then $F_i(\r)\sim F_{i+1}(\r)\#\r\p^2$. If 
$F_i\to F_{i+1}$ is the inverse of the blowing up
of a pair of conjugate points, then $F_i(\r)\sim F_{i+1}(\r)$. 
The induction works since $(S^1\times S^1)\#\r\p^2\sim \#3\r\p^2$.

Thus  we are reduced to
one of the following two cases:
\begin{enumerate}
\item   $F$ has a conic bundle structure $F\to B$, or
\item   $F$ is  Del Pezzo and $\rho(F)=1$.
\end{enumerate}

In the first case, $B_{\c}\cong \c\p^1$ since $F_{\c}$ is rational. 
Thus either
$B(\r)=\emptyset$ and so $F(\r)=\emptyset$, or $B\cong \r\p^1$.
Thus $F(\r)$ is the torus or the Klein bottle if there are no singular
fibers and
  $F(\r)\sim \uplus mS^2$ if there are $2m>0$ singular fibers
by (\ref{mmp.over.R.top}.C). Note that $S^2=\#0\r\p^2$ by convention.

In the second case we use  (\ref{mmp.over.R.top}.D). \qed

\begin{exrc}\label{comp.bir.inv}
 Let $X$ and $Y$ be smooth projective varieties over $\r$.
Assume that $X$ and $Y$ are birational to each other (over $\r$). Show
that $X(\r)$ and $Y(\r)$ have the same number of connected components.
\end{exrc}

\begin{say}[Vector bundles over real varieties] Let $X$ be a veriety
over $\r$ and $p:V\to X$ a vector bundle of rank $n$. Locally
$V$ is like
$U\times {\mathbb A}_{\r}^n\to U$ where $U\subset X$ is Zariski open
and ${\mathbb A}_{\r}^n=\spec_{\r}\r[t_1,\dots,t_n]$ is affine $n$-space
over $\r$. (Which should {\em not} be identified with $\r^n$!)

As usual, to $V$ one can associate a vector bundle $p_{\c}:V_{\c}\to
X_{\c}$ and also a real vector bundle
$p(\r):V(\r)\to X(\r)$ which is obtained by taking the $\r$-valued
points of ${\mathbb A}_{\r}^n$ which is exactly $\r^n$.

(To complete the picture, any real vector bundle on a manifold can be
complexified, and $V(\r)\otimes_{\r}\c\cong V_{\c}|_{X(\r)}$.)
\end{say}

\begin{say}[Degrees of line bundles over $\r$ and $\c$]
\label{dolb.exrc}{\ }

Let $B$ be a smooth projective curve over $\c$  and $L$ a line
bundle on $B$.  Let $s$ be a  nonzero meromorphic section of
$L$.  The number of zeros minus the number of poles of $s$ on $B$
(counted with multiplicity) is called the {\it degree of $L$}.  Let $Y$
be a smooth projective variety over
$\c$ and
$L$ a line bundle on $Y$.  For any curve $B\subset Y$
the degree of $L|_B$ is defined. It is
 also
called the {\it intersection number} of $B$ and $L$ and denoted by
$(B\cdot L)$. 

Let $A\sim S^1$ be a compact 1--dimensional manifold  and $L$ a real
line bundle on $M$.  Let $s$ be a  nonzero  section of
$L$.  The number of zeros  of $s$ on $A$ 
(counted with multiplicity)  makes sense only mod 2.  If  $M$ is a
compact manifold
  and
$L$ a real line bundle on $M$ then for any 1-cycle $A\subset M$
we obtain the $\z_2$-valued  {\it intersection number} of $A$ and $L$.
It is denoted by
$(A\cap L)$. (To be precise, I should write
$(A\cap w_1(L))$ where $w_1(L)$ stands for the first
Stiefel--Whitney class of  $L$. This is a class
in $H^1(X(\r),\z_2)$  analogous to the first Chern class of a
complex line bundle, cf. \cite[Sec. 4]{milnor-s74}.)

Let now $X$ be a smooth projective variety over $\r$, $C\subset X$ a
curve and $L$ a line bundle on $X$.  
We obtain two numbers:
$$
(L(\r)\cap C(\r)) \qtq{and}  (C_{\c}\cdot L_{\c}).
$$
What is the relationship between them?

To answer this, take a real meromorphic section 
$s$ of $L$  which has only finitely many zeros and poles on $C$. 
When we count the real zeros  and poles of $s$ on $C(\r)$, we miss the
complex zeros  and poles of $s$ on $C_{\c}$. Since $s$ is real, the
complex zeros  and poles come in conjugate pairs. Thus we conclude that
$$
(C(\r)\cap L(\r))\equiv (C_{\c}\cdot L_{\c}) \mod 2,
$$
which is best possible  since the left hand side is defined only mod 2
anyhow. 
\end{say}

\begin{say}[Orientability of $X(\r)$ and the canonical class]{\ }

Let $M$ be a differentiable manifold,
 $0\in M$  a point and $x_1,\dots,x_n$ local coordinates.
A {\it local orientation} of $M$ at $0$ is a choice of an $n$-form
$f(x)dx_1\wedge \dots\wedge dx_n$ with $f(0)\neq 0$ up to multiplication
by a positive function. An {\it orientation} of $M$ is a nowhere zero
global $n$-form  on $M$, up to multiplication
by a positive function. $n$-forms are sections of the real line bundle
$\det T^*_M$. 
If $S^1\sim A\subset M$ is a loop then one can choose a consistent
oreintation of $M$ along $A$ $\Leftrightarrow$  $\det T^*_M$ has a
nowhere zero section along $A$ $\Leftrightarrow$ $(\det T^*_M\cap A)=0$.

If $X$ is a smooth variety over $\r$ then $n$-forms appear as
sections of the canonical line bundle. This proves that
$$
\det T^*_{X(\r)} \cong K_X(\r).
$$
In many cases this gives a way to decide if $X(\r)$ is orientable or
not.
\end{say}

\begin{cor} Let $X$ be a smooth projective variety over $\r$. Assume
that there is a curve $C\subset X$ such that $(C\cdot K_X)$ is odd.
Then $X(\r)$ is not orientable.
\end{cor}

Proof. We have proved above that
$$
(\det T^*_{X(\r)}\cap C(\r))\equiv (C\cdot K_X) \equiv 1\mod 2.
$$
$C(\r)$ may have several components, but along one of them
$\det T^*_{X(\r)}$ has odd degree, so we can not choose a consistent
orientation along that component.\qed

\begin{exrc} Show that $\r\p^n$ is orientable iff $n$ is odd.
Let $X\subset \p^n$ be a smooth hypersurface of degree $d$.
Show that $X(\r)$ is orientable  if $n-d$ is odd.
Show that $X(\r)$ is not orientable  if $n$ and $d$ are both odd.
If $n$ and $d$ are both even, then $X(\r)$ may or may not be
orientable.
\end{exrc}

If $(C\cdot K_X)$ is even, then it can happen that $X(\r)$ is not
orientable along an even number of components of $C(\r)$. In some cases
we are still able to conclude orientability of $X(\r)$ using 
stronger assumptions:

\begin{exrc}\label{gen.orient.ex}
 Let $X$ be a smooth projective variety over $\r$. Assume
that $K_X\cong L^{\otimes 2}$ for a real line bundle $L$. 
Show that  $X(\r)$ is   orientable.

More generally, assume that $K_X\cong \o_X(2D+D')$ where $D,D'$ are
divisors over $\r$ and $D'(\r)$ has codimension at least 2 in $X(\r)$.
(This is equivalent to assuming that every irreducible component of
$D'$ is geometrically reducible.) Show that  $X(\r)$ is   orientable.
\end{exrc}

\section{Birational classification}

 \begin{defn} Let $F$ be a smooth real algebraic surface. A surface
obtained from $F$ by blowing up $a$ real points and $b$ pairs of
conjugate complex points (possibly infinitely near) is denoted by
$(F,a,2b)$.

Given $F$ and $a,b$, the surfaces of the form $(F,a,2b)$
consitute a connected family if $F(\r)$ and $F_{\c}$ are both connected.
\end{defn}

\begin{lem}\label{elem.bir.lem}
 We have the following elementary birational equivalences
between   the minimal models in (\ref{mmp.over.R}). 
\begin{enumerate}
\item $(\p^2,2,0)\cong (Q^{2,2},1,0)$.
\item $(\p^2,0,2)\cong (Q^{3,1},1,0)$.
\item $(Q^{4,0},0,2)$ is isomorphic to the blow up of $Q^{3,0}\times
\p^1$ at a pair of conjugate points on the same section 
$Q^{3,0}\times P$, $P\in \r\p^1$.
\item  Any minimal conic bundle over a rational curve with $2$ singular
fibers is isomorphic to $(Q^{3,1},0,2)$.
\end{enumerate}
\end{lem}

Proof. In the first two cases we blow up the 2 points in $\p^2$ and
then contract the line through them to get a quadric.

$Q^{4,0}\cong Q^{3,0}\times Q^{3,0}$, let $\pi_1$ be the first
projection. The pencil of planes through the 2 points gives a map
$p:Q^{4,0}\map \p^1$. 
$$
(\pi_1,p):Q^{4,0}\map Q^{3,0}\times \p^1 \qtq{is birational}
$$
and becomes a morphism after blowing up the 2 points.

Finally  assume that  $F\to B$ is a minimal conic bundle over a rational
curve with $2$ singular fibers. By (\ref{mmp.over.R}.C),
$B(\r)\neq\emptyset$, thus $B\cong \p^1$. $F_{\c}$ is the blow up of a
minimal ruled  surface  $F''$ at 2 points. We can even assume that
$F''$ has a  section $E$ with negative selfintersection $(E^2)=-k$
and the two points  are not on $E$. If $k\geq 2$ then all other
sections of $F''$ have selfintersection at least 2, so  $E\subset F$ is
 the unique section   with negative selfintersection. Thus $E$ is
defined over $\r$ and $F\to B$ is not minimal. 

Thus $k=1$ and there is a unique section $E'\subset F''$ such that
${(E'}^2)=1$ and $E'$ passes through the two blown up points.
Let $\bar E\subset F$ be the birational transform. Then $E$ and $\bar
E$ have to be conjugate. Contracting them gives  the quadric
$Q^{3,1}$.\qed

\begin{lem}\label{minruled.lem}
 Let $F$ be a smooth projective surface over $\r$ such that
$F_{\c}$ is a minimal ruled surface over $\c\p^1$.  Then
 one of the
following holds:
\begin{enumerate}
\item $F_{\c}\cong \p^1\times \p^1$ (these cases were enumerated in
(\ref{quadric.lem})), or
\item $F\cong \proj_B(\o_B+\o_B(-r))$ is a minimal ruled surface over a
smooth real conic $B$ for some $r>0$.
\end{enumerate}
\end{lem}

Proof.  By the classification of minimal ruled surfaces,
either $F_{\c}\cong \p^1\times \p^1$, or $F_{\c}$
has a unique irreducible curve 
$E$ with negative selfintersection $-r$. $E$ and the ruling $g:F\to B$
are then defined over
$\r$ and $g_*\o_F(E)\cong \o_B+\o_B(-r)$. \qed
\medskip

As a corollary, we obtain the following birational classification
of real surfaces such that $F_{\c}$ is rational:

\begin{cor}\label{birclass.over.R}
 Let $F$ be a smooth  real  projective surface
such that $F_{\c}$ is rational. Then $F$ is birationally equivalent
over $\r$ to a surface in exactly one of the following classes:
\begin{enumerate}
\item[1.] $Q^{3,0}\times \p^1$. In this case  $F(\r)=\emptyset$.
\item[2.] $\p^2$. In this case  $F(\r)$ is connected.
\item[3$_m$.]  Minimal conic bundle with $2m$ ($m\geq 2$) singular
fibers. In this case
$F(\r)$ has  $m$ connected components.
\item[4.] Degree 2 minimal Del Pezzo surface.
\item[5.] Degree 1 minimal Del Pezzo surface.
\end{enumerate}
\end{cor}

Remark. In (\ref{deg2.equiv.thm}) and  (\ref{deg1.equiv.thm})
we prove that  $F(\r)$ has  4  (resp. 5) connected components if $F$ is
a minimal Del Pezzo surface of degree 2 (resp. 1). 
\medskip

Proof. Let $F\to F^*$ be the minimal model of $F$. By (\ref{mmp.over.R})
$F^*$ is either one of those listed above, or $F^*$ is a conic bundle
with $0$ or $2$ singular fibers. 
The former are treated in (\ref{minruled.lem}).
The latter are birational to 
$Q^{3,0}\times \p^1$  by (\ref{elem.bir.lem}).

The number of connected components of the real part is a birational
invariant (\ref{comp.bir.inv}), hence the cases (1--3$_m$) are all
different birationally.

 The cases (4--5)  differ
birationally from the other ones by (\ref{seg-man.deg123}).\qed

We use, without proof, the following result about the birational
classification of low degree Del Pezzo surfaces over any field.
Lectures  2--3 of \cite{KS97}  serve as a good introduction.

\begin{thm}\cite{Segre51, Manin66}\label{seg-man.deg123}
Let  $k$ be a  field (of characteristic zero) and $F$ a   minimal Del
Pezzo surface of degree 1,2 or 3 over $k$. Then 
\begin{enumerate}
\item $F$  is not rational (over $k$),
\item $F$  is not birational (over $k$) to any conic fibration,
\item $F$  is  birational to another minimal Del
Pezzo surface  $F'$ of degree 1,2 or 3 over $k$ iff $F$ is isomorphic to
$F'$. \qed
\end{enumerate}
\end{thm} 

This theorem, (\ref{birclass.over.R}) and (\ref{comp.bir.inv}) imply the
following:

\begin{cor}\label{rat.char.over.R}\cite[VI.6.5]{Silhol89}
 Let $F$ be a smooth projective surface  over $\r$.
The following are equivalent:
\begin{enumerate}
\item $F$ is birational to $\p^2$ over $\r$.
\item  $F_{\c}$ is birational to $\c\p^2$ and $F(\r)$ is connected. \qed
\end{enumerate}
\end{cor}

\section{Birational Classification of Conic Fibrations}

\begin{defn} Let $F$ be a smooth projective surface over a field $k$. A
morphism
$f:F\to B$  to a smooth curve is called a {\it conic fibration} if the
generic fiber is isomorphic to a plane conic (over $k(B)$). 
By (\ref{conics.r}) this is equivalent to assuming that $f^{-1}(b)\cong
\p^1_{\bar k}$ for a general $b\in B(\bar k)$.
\end{defn}

In this section we discuss the birational classification
of those surfaces over $\r$ which admit a
 conic fibration. This covers all surfaces where the MMP ends with  the
case (\ref{mmp.surf.thm}.C).

This is done in two steps.  First we consider those
 birational maps which  preserve the conic fibration. To be precise:

\begin{defn} Two conic fibrations $f:F\to B$ and $f':F'\to B'$ are
called {\it birational} if there is a birational map
$\phi:F\map F'$   and an isomorphism $\tau: B  \cong  B'$
 (both over $k$) which give a  commutative diagram
$$
\begin{array}{ccc}
F & \stackrel{\phi}{\map} & F'\\
f\downarrow{\ } & & {\ }\downarrow f'\\
B & \stackrel{\tau}{\map} & B'
\end{array}
$$
\end{defn}

The second step is to understand the birational maps between $F$ and
$F'$ which do not preserve the conic fibration. Fortunately, in many
cases there are no such maps. (For a proof see
\cite[V--VI]{Silhol89} or the original paper of
\cite{Iskovskikh67}.) 

\begin{thm}\label{cb.unique.thm}
 Let $k$ be a field and $f:F\to B$ a relatively minimal
conic bundle over $B$. Let $f':F'\to B'$ be
any   conic fibration and  $\phi:F\map F'$  any birational map
 (over $k$). Then 
\begin{enumerate}
\item If  $(K_F^2)\leq 0$ then $\phi$ is a birational map of the
conic fibrations.
\item If  $(K_F^2)=2$ then $F$ and $F'$ are   birational 
conic fibrations (though $\phi$ itself need not respect the fibration
structure).\qed
\end{enumerate}
\end{thm}

\begin{defn} Let $f:F\to B$ be a conic fibration over $\r$.
The image of the set of real points 
$f(F(\r))\subset B(\r)$ is   a union of finitely many closed intervals. 
Let us denote it by $I(F)$.
\end{defn}

The main theorem of the section shows that $I(f)$ characterizes $f$:

\begin{thm} Two conic fibrations $f:F\to \p^1$ and $f':F'\to \p^1$
over $\r$ are
 birational iff there is  an isomorphism $\tau:\p^1  \cong  \p^1$
such that $\tau(I(f))=I(f')$. 
\end{thm}

Proof. Let $\phi:F\map F'$ and $\tau:\p^1  \cong  \p^1$ be a
birational map of the two 
 conic fibrations. Then
$F'(\r)$ and $\phi(F(\r))$ agree outside finitely many fibers, thus 
$I(f')$ and $\tau(I(f))$   differ only at finitely many points. 
Unions of closed intervals can not differ at finitely many points only,
thus in fact $\tau(I(f))=I(f')$.

The converse is established by bringing each conic fibration to a
normal form. (The roots $a_i$ in (\ref{cb.standard.form.thm})
are the boundary points of $I(f)$. This leaves
two choices for $I(f)$ itself, corresponding to the two choices of
the sign on the right hand side.)
\qed

\begin{thm}\label{cb.standard.form.thm}
 Let $f:F\to \p^1$ be a conic fibration over $\r$. Then
$f$ is birational   to a conic fibration $f':F'\to \p^1$ with affine
equation
$$
x^2+y^2=\pm \prod_{i=1}^{2m}(z-a_i)\subset {\mathbb A}^3,
$$
where the $a_i$ are distinct real numbers.
\end{thm}

The proof rests on the following simple lemma about quadratic forms:

\begin{lem}\label{cb.standard.form.lem} Let $k$ be a field (of
characteristic different from 2) and $Q(x_0,\dots,x_n)$ a quadratic form
over $k$ which is anisotropic (that is $Q=0$ has no nontrivial solution
over $k$). For any $a\in k$ the following are equivalent
\begin{enumerate}
\item $Q=0$ has a nontrivial solution over $k(\sqrt{a})$.
\item After a suitable coordinate change, $Q$ can be written as

\noindent $b(y_0^2-ay_1^2)+Q'(y_2,\dots,y_n)$. 
\end{enumerate}
\end{lem}

Proof. (2) $\Rightarrow$ (1) is shown by the substitution
${\bf y}:=(\sqrt{a}:1:0\cdots:0)$. 

Conversely, assume that ${\bf v}\in k(\sqrt{a})^{n+1}$ 
satisfies $Q({\bf v})=0$. Let $\bar {\bf v}$ denote the
conjugate of ${\bf v}$. Then $Q(\bar{\bf v})=0$.
${\bf v}$ and $\bar{\bf v}$ span a 2-dimensional linear
subspace of $k(\sqrt{a})^{n+1}$  which is defined over $k$. That is,
there is a linear subspace $V\subset k^{n+1}$ such that
$V\otimes_kk(\sqrt{a})=\langle {\bf v}, \bar{\bf
v}\rangle$.  $Q$ is nondegenerate on $V$ (since $Q=0$ has no solutions
in $k$), thus $k^{n+1}=V+V^{\perp}$ where $V^{\perp}$ is the orthogonal
complement of $V$ with respect to $Q$. Let $y_2,\dots,y_n$ be
coordinates on $V^{\perp}$ and choose coordinates $y_0,y_1$ on $V$
such that 
$$
{\bf v}=(\sqrt{a},1)\qtq{and} \bar{\bf v}=(-\sqrt{a},1).\qed
$$

\begin{say}[Proof of (\ref{cb.standard.form.thm})]{\ }

Let $k:=\r(t)$ be the quotient field of $\p^1_{\r}$. The generic fiber
of $F\to \p^1$ is birational to  a plane conic $C_k$ over $k$
(\ref{conics.r}).

If $C_k$ has a $k$-point (equivalently, if $F\to \p^1$ has a section)
then $F$ is birational to $\p^1\times \p^1$ by (\ref{conics.r}).  $C_k$
has a point over
$k(\sqrt{-1})=\c(t)$ (equivalently,  $F\to \p^1$ has a section over
$\c$). Thus by  (\ref{cb.standard.form.lem}), in suitable coordinates
the equation of $C_k$ becomes $x_0^2+x_1^2=g(t)x_2^2$
for some $g(t)\in \r(t)$. 

We can multiply through with the square of the denominator of $g$,
thus we may assume that $g(t)\in \r[t]$. 

Write $g=f(t)g_1(t)^2\prod_a (z-a)(z-\bar a)$ where $f$ has only simple
real roots and the $a$ are nonreal complex numbers. We can divide by 
$g_1(t)^2$. If $a=u+iv$ then $(z-a)(z-\bar a)=(z-u)^2+v^2$. 
Note that
\begin{eqnarray*}
g(t)(h_0^2+h_1^2)x_2^2&=&x_0^2+x_1^2\qtq{is equivalent to}\\
g(t)x_2^2&=&
\left(\frac{x_0h_0-x_1h_1}{h_0^2+h_1^2}\right)^2+
\left(\frac{x_0h_1+x_1h_0}{h_0^2+h_1^2}\right)^2.
\end{eqnarray*}
Using this, we can get rid of the complex factors $\prod_a (z-a)(z-\bar
a)$ one at a time. At the end we obtain the required normal form,
except that we may have an odd number of factors on the right hand
side:
$$
x^2+y^2=\pm \prod_{i=1}^{2m-1}(z-a_i).
$$
 In this case we first apply a
translation to ensure that $0$ is not among the $a_i$ and then make a
substitution
$(x,y,z)\mapsto (xz^{-n},yz^{-n}, z^{-1})$ to get the equation
$$
x^2+y^2=\pm \prod_{i=1}^{2m}(z-a'_i),
$$
where $a'_i=a_i^{-1}$  for $i<2m$ and  $a'_{2m}=0$.\qed
\end{say}

Putting things together,  we obtain the following criterion for
birational equivalence of conic   fibrations.
The result corrects a slight inaccuracy in \cite[VI.3.15]{Silhol89}.

\begin{exrc}\label{cb.isom.ex}  Two conic bundles
\begin{eqnarray*}
F & = & (x^2+y^2= c\prod_{i=1}^{2m}(z-a_i))\to \p^1\qtq{and}\\
F' & = & (x^2+y^2= c'\prod_{i=1}^{2m}(z-a'_i))\to \p^1
\end{eqnarray*}
are birational to each other iff there is  a permutation $\sigma\in
S_{2m}$  and  a matrix
$$
\left(\begin{array}{cc}
\alpha & \beta \\
\gamma & \delta
\end{array}
\right)\in GL(2,\r)
\qtq{such that}
a'_{\sigma(i)}=\frac{\alpha a_i+\beta }{\gamma a_i+\delta}, \qtq{and}
$$
$c'$ and $c\prod_i(\gamma a_i+\delta)$ have the same sign.
\end{exrc}

\begin{exrc}\label{geom.cb.ex}
 Using elementary transformations of conic bundles, give a
geometric proof of the results in this section.
\end{exrc}

\begin{say}[Moduli of conic fibrations]{\ }

Let $F$ be a smooth projective surface over $\r$ which admits a conic
fibration $f:F\to \p^1$.  We   proved that if $I(f)\subset \r\p^1$
has at least 3 components then $I(f)$ (modulo the action of $GL(2,\r)$)
determines $F$ up to birational equivalence (over $\r$). 

The space of $m$ disjoint closed intervals in $\r\p^1$ is a connected
manifold of real dimension $2m$. The quotient by the $GL(2,\r)$ action
gives a $2m-3$ dimensional topological space (it has some quotient
singularities).  With some more care, we could even realize this space
as the set of real points of a $(2m-3)$-dimensional algebraic variety.

For $m=0$ there are 2 conic fibrations up to birational equivalence:
$\p^1\times \p^1$ gives $I(f)=\r\p^1$ and $Q^{3,0}\times \p^1$ gives
$I(f)=\emptyset$. For $m=1$ we have only one birational equivalence
class by (\ref{elem.bir.lem}.4). 

For $m=2$ we see in (\ref{deg4.cb.ex}) that  all such surfaces are
birational to each other (though they are not birational as conic
fibrations). 
\end{say}

\section{Del Pezzo Surfaces of Degree $\geq 3$}

In this section we describe all Del Pezzo surfaces of degree $d\geq 3$
over
$\r$.

\begin{exrc}\label{amp.image.ex}
 Let $g:S\to S'$ be a birational morphism of smooth
surfaces. Show that if $H$ is ample on $S$ then $f(H)$ is ample on
$S'$. Thus if $S$ is Del Pezzo then $S'$ is also Del Pezzo. 
\end{exrc}

\begin{prop}\label{bireg.empty}
 Let $F$ be a smooth  real  Del Pezzo surface which is
birational to $Q^{3,0}\times \p^1$. 
Set $d:=(K_F^2)$. Then $d\in\{8,6,4,2\}$.

If $d=8$ then 
 $F$ is isomorphic to
either $Q^{3,0}\times Q^{3,0}$ or to $Q^{3,0}\times \p^1$. 

If $d<8$ then 
 $F$ is isomorphic to
 to $Q^{3,0}\times \p^1$  blown up in $\frac12(8-d)$ pairs of conjugate
points.

Therefore, for $d\in\{6,4,2\}$  such surfaces form a connected family.
\end{prop}

Proof.  Apply the MMP over $\r$ to obtain  $F\to \cdots\to F^*$. 
If $F$ is Del Pezzo then so is $F^*$ by (\ref{amp.image.ex}),
and $(K_{F^*}^2)\geq (K_F^2)$. Hence
in our case $F^*$  is either 
$Q^{3,0}\times Q^{3,0}$ or  $Q^{3,0}\times \p^1$. 
By (\ref{elem.bir.lem}.3) any blow up of $Q^{3,0}\times Q^{3,0}$ at 
a pair of conjugate points is also a blow up of 
$Q^{3,0}\times \p^1$. \qed

\begin{prop}\label{bireg.rtl} Let $F$ be a smooth  real  Del Pezzo
surface which is birational to $\p^2$. 
Set $d:=(K_F^2)$. Then $9\geq d\geq 1$ and we have one of the following
cases:

If $d=9$ then 
 $F$ is isomorphic to $\p^2$.

If $d=8$ then 
 $F$ is isomorphic to
either $Q^{3,1}$ or to $Q^{2,2}$ or to 
$\p^2$  blown up at a real point.

If $d<8$ then 
 $F$ is isomorphic to one of the following:
\begin{enumerate}
\item 
 $\p^2$  blown up at $a\geq 0$ real points and $b\geq 0$ pairs of
conjugate points for some $a+2b=9-d$. Thus $F(\r)\sim \#(a+1)\r\p^2$.
\item $Q^{3,1}$  blown up at  $b=\frac12(8-d)$ pairs of
conjugate points (so $d$ is even). Thus $F(\r)\sim  S^2$.
\item $Q^{2,2}$  blown up at  $b=\frac12(8-d)$ pairs of
conjugate points (so $d$ is even). Thus $F(\r)\sim  S^1\times S^1$.

Therefore, for any $d<8$, such surfaces with a given topological
type $F(\r)$ form a connected family.
\end{enumerate}
\end{prop}

Proof. The minimal model of such a surface is either  $\p^2$,
$Q^{3,1}$ or  $Q^{2,2}$. 
By (\ref{elem.bir.lem}.1--2) any blow up of $Q^{3,1}$ or to $Q^{2,2}$ at 
a real point is also a blow up of 
$\p^2$. \qed
\medskip

The two propositions above account for all Del Pezzo surfaces of degrees
$d\geq 5$. The results are summarized in the next statement:

\begin{cor}\label{d>4.dp.top} The following table lists all topological
types of the real points of Del Pezzo surfaces of degrees $9\geq d\geq
5$. All surfaces of a fixed degree and topological type form a
connected family, except for $d=8$ and $F(\r)=\emptyset$ when there are
2 such surfaces.
$$
\begin{tabular}{|c|l|}
\hline degree & \qquad\qquad topological types\\
\hline 9 & $\r\p^2$\\
\hline 8 & $S^2$ or $S^1\times S^1$ or $\r\p^2\#\r\p^2$ or $\emptyset$ \\
\hline 7 & $\r\p^2$ or $\#3\r\p^2$ \\
\hline 6 & $S^2$ or $S^1\times S^1$ or $\r\p^2\#\r\p^2$  or
        $\#4\r\p^2$ or $\emptyset$\\
\hline 5 & $\r\p^2$ or $\#3\r\p^2$ or $\#5\r\p^2$ \\
\hline
\end{tabular}
$$ 
\end{cor}

The following result shows that odd degree Del Pezzo surfaces over $\r$
are relatively easy to understand:

\begin{lem}\label{d1=d2.blowup}
 Every  degree $2d-1$  Del Pezzo surface $F$ over $\r$  with
$\rho(F)\geq 2$ is the blow up of a degree $2d$ Del Pezzo surface at a
real point.
\end{lem}

Proof.  Since $(K_F^2)$ is odd, $F$ is not a minimal conic bundle.
Thus $F$ is either the  
blow up of a degree $2d$ Del Pezzo surface at a real point
or the blow up of a degree $2d+1$ Del Pezzo $F'$ surface at a conjugate
pair of complex points $P+\bar P$.  

If $F'\cong \p^2$ then let $L\subset \p^2$ be the line through the two
points. Its birational transform on $F$ is a line.

Otherwise, $F'$ is again not minimal by (\ref{mmp.over.R}), hence 
$F'$ contains a line $L$ over $\r$ by induction.
$P,\bar P\not\in L$ since otherwise the birational transform of $L$ on
$F$ would have a nonnegative intersection number with $K_F$.
Thus $L$ gives a real line on $F$. 

Contracting a real line on $F$   we get a degree $2d$ Del
Pezzo surface.\qed

\medskip

This shows that the study of degree 3   Del Pezzo surfaces
is reduced to the study of  degree 4 cases. 
The classification of these two classes is summarized next.
 These results were   obtained by
\cite{Schlafli1863}, who actually worked directly with cubic surfaces.

\begin{cor}\label{d=4,3.dp.top} The following table lists all topological
types of the real points of Del Pezzo surfaces of degrees 4 and 3. All
surfaces of a fixed degree and topological type form a connected family.
$$
\begin{tabular}{|c|l|}
\hline degree & \qquad\qquad topological types\\
\hline 4 & $S^2$ or $S^1\times S^1$ or $\r\p^2\#\r\p^2$  or
        $\#4\r\p^2$ or $\emptyset$ or $S^2\uplus S^2$\\
\hline 3 & $\r\p^2$ or $\#3\r\p^2$ or $\#5\r\p^2$ or $\#7\r\p^2$
or $S^2\uplus \r\p^2$ \\
\hline
\end{tabular}
$$ 
 Moreover, in the $S^2\uplus S^2$ case the
monodromy interchanges the two components.
\end{cor}

Proof.  As we noted above,  it is sufficient to describe
all degree 4 
Del Pezzo surfaces.

If a degree 4 
Del Pezzo surface $F$ is obtained from a  higher degree surface by
blowing up then we are reduced to (\ref{d>4.dp.top}). Otherwise
$F$ is a conic bundle over $\p^1$ with 4 singular fibers.
These 4 singular fibers give 8 lines on $F$. By looking at the set of
all lines over $\c$, we see that the remaining 8 lines again form 4
pairs and determine another morphism to $\p^1$. Thus $F$ is a double
cover of $\p^1\times \p^1$ ramified along a curve $D\subset \p^1\times
\p^1$ of type
$(2,2)$. $D$ has 4 horizontal and 4 vertical tangents and the 16 lines
are   sitting over these tangents.

The rest is a special case of (\ref{(2,2)-curves}). \qed

\begin{exrc}\label{(2,2)-curves}
Show that the space of all smooth real curves of type $(2,2)$ on
$\p^1\times
\p^1$ has 7 connected components. They are determined by the homotopy
classes of the components of $D(\r)$: $\emptyset$ or $(0,0)$ or
$(0,0)\uplus (0,0)$ or $(1,1)\uplus (1,1)$ or $(1,-1)\uplus (1,-1)$ or
$(1,0)\uplus (1,0)$ or $(0,1)\uplus (0,1)$.
\end{exrc}

\begin{exrc}\label{deg4.cb.ex}\cite[VI.3.5]{Silhol89} 
  Using the correspondence between $(2,2)$-curves 
on
$\p^1\times\p^1$ and degree 4 Del Pezzo surfaces 
show that any two minimal conic bundles $F,F'$  with 4 singular fibers
are birational, by producing examples  $S\to \p^1\times \p^1$ such that
one conic bundle structure  of $S$ is birational to $F$ and the other
to $F'$.
(This should be easier after the next section.)
\end{exrc}

\begin{exrc}\cite{Schlafli1863}\label{cubic.lines.ex}
  Show that a smooth cubic over $\r$ has 27,15,7 or 3 real
lines.
\end{exrc}

\section{Del Pezzo Surfaces of Degree 2 and 1}

\begin{notation} Let $D\subset \p^2$ be  a degree 4 smooth real curve.  
$D(\r)$ divides
$\r\p^2$ into connected open sets and precisely  one of these  is
nonorientable (denoted by $U_D$).  We choose an equation $f(x,y,z)\in
\r[x,y,z]$ of $D$  such that $f$ is negative
on $U_D$.

We can associate two different degree 2  Del Pezzo surfaces to $D$. 
One is $F^+_D:=(u^2=f(x,y,z)\subset \p^3(1,1,1,2)$
and the other $F^-_D:=(u^2=-f(x,y,z)\subset \p^3(1,1,1,2)$.

The correspondence $F^+_D\leftrightarrow F^-_D$
is a natural involution on the space of degree 2 real Del Pezzo
surfaces.

$D$ has 28 bitangents over $\c$ and over each bitangent of $D$
we get a pair of lines on $F^{\pm}_D$. This gives a total of 56 lines.
\end{notation}

The topological classification of degree 4 plane curves over $\r$ is
very old, it is already contained in \cite{Plucker1839}.
(See \cite{Viro90} for a recent survey of the study of low degree
real plane curves.)
This implies the topological classification of degree 
2 real Del Pezzo surfaces. The following proposition summarizes
these results.

\begin{prop}\label{deg2.dp.top}
 There are 6 topological types of degree 4 smooth real plane
curves. Correspondingly
there are 12 topological types of degree 2 real Del Pezzo surfaces.
The following table gives the complete list.
The types
in the same row correspond to each other under 
$D \leftrightarrow F^+_D\leftrightarrow F^-_D$.
$$
\begin{tabular}{|c|c|c|}
\hline $D(\r)$ & $F^+_D(\r)$  & $F^-_D(\r)$\\
\hline $\bigcirc\bigcirc\bigcirc\bigcirc{}$ & $\uplus 4S^2$ &
             $\#8\r\p^2$\\
\hline $\bigcirc\bigcirc\bigcirc$ & $\uplus 3S^2$ &  $\#6\r\p^2$\\
\hline $\bigcirc\bigcirc{}$ & $S^2\uplus S^2$ &  $\#4\r\p^2$\\
\hline $\bigcirc$ & $S^2$ &  $ \#2\r\p^2$\\
\hline $\emptyset$ & $\emptyset$ &  $\r\p^2\uplus \r\p^2$\\
\hline $\bigcirc \!\!\!\!\!\circ$ & $S^1\times S^1$ & $S^2\uplus
 \#2\r\p^2$\\
\hline
\end{tabular}
$$ 
\end{prop}

\cite{Zeuthen1874} studied the bitangents of degree 4 plane curves.
He proved the equivalence of (\ref{deg2.equiv.thm}.1) and
(\ref{deg2.equiv.thm}.5).  He understood the relationship between
degree 4 plane curves and cubic surfaces. 
(Projecting a cubic surface from one of its points, the branch curve
is a plane quartic. Equivalently, blowing up the cubic at a point we
get a degree 2 Del Pezzo surface.)
This is, however, not the
 natural thing to do from the modern viewpont.
Most of (6.3) is proved in \cite{Comessatti13}.

\begin{thm} \label{deg2.equiv.thm}
Let $D\subset \p^2$ be  a degree 4 smooth real curve.  The
following are equivalent:
\begin{enumerate}
\item All 28 bitangents of $D$ are real.
\item All 56 lines of $F^-_D$ are real.
\item $F^-_D$ is isomorphic to $\p^2$ blown up in 7 real points.
\item $F^-_D(\r)\sim \#8\r\p^2$.
\item $D(\r)\sim \uplus 4S^1$.
\item  $F^+_D(\r)\sim \uplus 4S^2$.
\item  $F^+_D$ has Picard number 1 over $\r$.
\end{enumerate}
\end{thm}

Proof. (1) $\Rightarrow$ (2): A neighborhood of a line in $\r\p^2$ is
not orientable, thus any bitangent is contained in $U_D$ (except for the
points of tangency).   $f$ is negative on any bitangent and so
$u^2=-f$ has real solutions, giving 56 real lines on 
$F^-_D$.

(2) $\Rightarrow$ (3): Over $\c$, $F^-_D$ is the blow up of $\p^2$ at 7
points, hence it has 7 disjoint lines. If all lines are real, we have 7
disjont real lines. Contracting these we get a Del Pezzo surface of
degree 9 over
$\r$. By (\ref{mmp.over.R}.D) it is $\p^2_{\r}$. 

(3) $\Rightarrow$ (4): Topologically, each blowing  up is connected sum
with
$\r\p^2$.

(4) $\Rightarrow$ (5): This follows from (\ref{deg2.dp.top}).

(5) $\Rightarrow$ (6): This also follows from (\ref{deg2.dp.top}).

(6) $\Rightarrow$ (7): Assume to the contrary that 
$F^+_D$ has Picard number $\geq 2$ over $\r$. By (\ref{mmp.over.R})
we have one of 2 cases:
\begin{enumerate}
\item $F^+_D$ is a minimal conic bundle with 6 singular fibers. In this
case $F^+_D(\r)\sim \uplus 3S^2$, a contradiction.
\item  $F^+_D$ is the blow up of a   Del Pezzo surface of degree 3
or 4 over
$\r$. By (\ref{d=4,3.dp.top}) 
$F^+_D(\r)$ has at most 2 connected components, a contradiction.
\end{enumerate}

(7) $\Rightarrow$ (1):  Assume that $D$ has a complex bitangent $L$.
Its conjugate $\bar L$ is again a bitangent. Let $C\subset F^+_D$ be a
complex line over $L$. Its conjugate $\bar C$ lies over $\bar L$.
Then $(C\cdot \bar C)\leq 1$ (the only possible intersection point lies
over $L\cap \bar L$). Thus $F^+_D$ has either a disjoint pair of
conjugate lines or a conic bundle structure, a contradiction.\qed

\begin{prop}\label{deg2.dp.top.conn}\cite{Klein1876}
 The space of degree 4 smooth real plane
curves has 6 connected components corresponding to the
6 topological types in (\ref{deg2.dp.top}). The space  of degree 2
real Del Pezzo surfaces has 12 connected components corresponding to the
12 topological types in (\ref{deg2.dp.top}).
\end{prop}

Proof. The two parts are equivalent and it is sufficient to treat the
$F^-_D$ cases.  $D$ always has  a real bitangent (\ref{deg.4.bitang}),
thus 
$F^-_D$ contains a real line. So $F^-_D$ is obtained by blowing up a
degree 3 Del Pezzo surface at a real point. 
(\ref{deg2.dp.top.conn}) now follows from (\ref{d=4,3.dp.top}).\qed

\begin{exrc}\cite{Zeuthen1874}\label{deg.4.bitang}
 Let $C$ be  a smooth real plane curve
of degree 4. Let $d$ be the number of ovals of $D(\r)$
which are not contained in another oval.
Show that $C$ has $4+2d(d-1)$ real bitangents.
\end{exrc}

\begin{say}[Degree 1 Del Pezzo surfaces]\label{deg1.dp.say}{\ }

Let $F$ be a degree 1 Del Pezzo surface over any field $k$.  $|-K_F|$ is
a pencil with exactly one base point. So this is a $k$-point  and
$F(k)\neq
\emptyset$. $|-2K_F|$ is base point free and exhibits $F$ as a double
cover of a quadric cone $Q\subset \p^3$, ramified along a  curve
$D\subset Q$ which is a complete
intersection of $Q$ with a cubic surface with equation $(f=0)$. 
$D$ does not pass through the vertex of the cone. 

$F_{\bar k}$ contains 240 
lines (that is $-1$-curves);  cf. \cite[IV.4.3]{Manin72}. We obtain these
  as follows. Take a  plane $H\subset \p^3$ which is tangent
to
$D$ at 3 points. The preimage of $H\cap Q$ in $F$ has 2 irreducible
components, each is a line. Thus we conclude that there are 120 planes 
which are tangent to $D$ at 3
points.

Assume now that $k=\r$. 
Since $Q(\r)\neq \emptyset$, we can write 
$Q=(x^2+y^2=1)$ 
in suitable affine coodinates 
$(x,y,z)$ on ${\Bbb A}^3$. That is, $Q(\r)$ is a cylinder with a singular
point at infinity.

As in the degree 2 case, for each (nonhomogeneous) cubic $f(x,y,z)$ we
obtain two degree 1 Del Pezzo surfaces, given by  affine  equations
$$
F^{\pm}_f:=(x^2+y^2-1=u^2\mp f(x,y,z)=0)\subset {\Bbb A}^4.
$$

There are 2 types of  simple closed loops on a cylinder: null
homotopic ones (I  call them ovals) and those homotopic to a plane
section (I call them big circles).

Since $D(\r)$ is the intersection of a cylinder with a cubic, it has 3
or 1 intersection points with any ruling line of the cylinder.
Thus $D(\r)$ contains either  3 big circles (and no ovals) or 1 big
circle. $D$ has genus 4, hence by Harnack's theorem
(cf. \cite[VII.4]{Shafarevich72}), $D(\r)$ has at most 5 connected
components. An oval can not be inside another oval since this would give
4 points on a ruling. Furthermore, we can not have an oval on one side
the big circle and at least two ovals on the other side. Indeed, choosing
points
$P_1,P_2,P_3$ inside the 3 ovals, there is a plane $H$ through them. 
Then $H$ intersects each oval in at least 2 points, and also the big
circle. So $(H\cdot D)\geq 8$, but $D$ has degree 6, a contradiction.

If all the ovals are on the same side of the big circle, we can
normalize $f$ so that it is positive on the other side.
The other cases are symmetrical and it makes little sense to normalize
$f$.

We can summarize these results in the following table:
\end{say}

\begin{prop}\label{deg1.dp.top}
 There are 7 topological types of degree 6 smooth real complete
intersection curves on the cylinder  $(x^2+y^2=1)$, not passing through
the point at infinity. Correspondingly there are 11 topological types of
degree 1 real Del Pezzo surfaces. The following table gives the complete
list. The types
in the same row correspond to each other under 
$D=(f=0)\cap Q \leftrightarrow F^+_f\leftrightarrow F^-_f$.
$$
\begin{tabular}{|c|c|c|}
\hline $D(\r)$ & $F^+_f(\r)$  & $F^-_f(\r)$\\
\hline 1 big circle + 4 ovals & $\r\p^2\uplus 4S^2$ & 
           $ \#9\r\p^2$\\
\hline 1 big circle + 3 ovals & $\r\p^2\uplus 3S^2$ & 
           $ \#7\r\p^2$\\
\hline 1 big circle + 2 ovals & $\r\p^2\uplus 2S^2$ & 
           $ \#5\r\p^2$\\
\hline 1 big circle + 1 oval & $\r\p^2\uplus  S^2$ & 
           $ \#3\r\p^2$\\
\hline 1 big circle + 0 oval & $\r\p^2$ &  $\r\p^2$\\
\hline 1 big circle + 1+1 ovals & $ \#3\r\p^2\uplus  S^2$ & 
           $ \#3\r\p^2\uplus  S^2$\\
\hline 3 big circles  & $\r\p^2\uplus  \#2\r\p^2$ &
             $\r\p^2\uplus  \#2\r\p^2$\\
\hline
\end{tabular}
$$ 
\end{prop}

The following theorem, due to  \cite{Comessatti13}, is the degree 1
version of (\ref{deg2.equiv.thm}).  I thank F. Russo for checking the
arguments of Comessatti.

\begin{thm} \label{deg1.equiv.thm}
Let $D=(f=0)\subset Q$ be  a degree 6 smooth real complete
intersection curve on the cylinder  $Q=(x^2+y^2=1)$.  The following are
equivalent:
\begin{enumerate}
\item All 120 triple tangents of $D$ are real and $f$ is negative on
all of them.
\item All 240 lines of $F^-_f$ are real.
\item $F^-_f$ is isomorphic to $\p^2_{\r}$ blown up in 8 real points.
\item $F^-_f(\r)\sim \#9\r\p^2$.
\item $D(\r)\sim \uplus 5S^1$.
\item  $F^+_f(\r)\sim \r\p^2\uplus 4S^2$.
\item  $F^+_f$ has Picard number 1 over $\r$.
\end{enumerate}
\end{thm}

Proof. (1) $\Rightarrow$ (2):  If $f$ is negative on a triple tangent
then
$u^2=-f$ has real solutions, giving a pair of real lines on 
$F^-_f$.

(2) $\Rightarrow$ (3): Over $\c$, $F^-_f$ is the blow up of $\p^2$ at 8
points. Thus it has 8 disjoint lines. If all lines are real, we have 8
disjont real lines. Contracting these we get a Del Pezzo surface of
degree 9 over
$\r$. By (\ref{mmp.over.R}.D) it is $\p^2_{\r}$. 

(3) $\Rightarrow$ (4): Topologically, each blowing  up is connected sum
with
$\r\p^2$.

(4) $\Rightarrow$ (5): This follows from (\ref{deg1.dp.top}).

(5) $\Rightarrow$ (6): This also follows from (\ref{deg1.dp.top}).

(6) $\Rightarrow$ (7): Assume to the contrary that 
$F^+_f$ has Picard number $\geq 2$ over $\r$. 
$F^+_f$ can not be a minimal conic bundle since $(K^2)$ is odd.
Thus  $F^+_f$ is the blow up of   Del Pezzo surface of degree 2 or 3
over
$\r$. By (\ref{deg2.dp.top}, \ref{d=4,3.dp.top}) 
$F^+_f(\r)$ has at most 4 connected components, a contradiction.

(7) $\Rightarrow$ (1):  If $D$ has a complex triple tangent, we can
argue as in (\ref{deg2.equiv.thm}.(7) $\Rightarrow$ (1)) that $F^+_f$
contains a conjugate pair of lines $C,\bar C$ such that $(C\cdot \bar
C)\leq 2$.  $(C+\bar C)\equiv rK$ for some $r\in \z$, thus
$2(C\cdot \bar C)-2=r^2$. This is impossible.

If there is a real triple tangent such that  $f$ is positive on it then
as in (1) $\Rightarrow$ (2) we get real lines on  $F^+_f$.
\qed

\begin{prop}\label{deg1.dp.top.conn}
 The space  of degree 1
real Del Pezzo surfaces has 11 connected components corresponding to the
11 topological types in (\ref{deg1.dp.top}).
\end{prop}

Proof. For the  nonminimal ones, this follows from 
(\ref{deg1.dp.top.conn}) and  (\ref{d1=d2.blowup}). 
The minimal ones are in one--to-one correspondence with the blow ups of
$\p^2$ at 8 real points, this is again a connected space.\qed

\begin{ack} I thank my audience at the Trento sumer school, especially
L. Bonavero,  S. Cynk and  S. Endrass  for numerous comments and
improvements.   A. Marin directed me to several 19th century references.
F. Russo checked   the arguments of Comessatti  about degree 1 and 2
Del Pezzo surfaces and pointed out some misunderstandings on my part.
  Partial financial support was provided by  the NSF
under grant number  DMS-9622394. 
\end{ack}

\section{Hints to selected exercises}

(\ref{er.noncontr.ex}) Blow up $\geq 10$ general points (over $\c$) on a
smooth plane cubic. The birational transform of the cubic generates an
extremal ray which is not contractible. To see this show that $\pic$ of
the blown up surface injects into   $\pic$ of the cubic.
\medskip

(\ref{int.cone.exrc})  Write down Riemann--Roch for $mz$ and $-mz$ to
get the first part. Then use this for $z+\epsilon z'$ for any $z'\in
N_1(F)$. 
(cf. \cite[V.1.8]{Hartshorne77}.)
\medskip

(\ref{conics.r}) $|-K|$ embeds $C$ as a conic.
\medskip

(\ref{projspace.ex})  Let $H\in X_{\c}$ be a hyperplane. Show that
$H\cap \bar H$ is real and use induction.
Even degree symmetric powers of the empty conic give examples. 
\medskip

(\ref{coh.b.c}) This is a special case of 
\cite[III.9.3]{Hartshorne77}.
\medskip

(\ref{pic=gal.inv})  The empty conic gives a good example. 

For the rest, the key point is to understand that we know more than the
existence of an isomorphism $\tau:L_{\c}\cong \bar L_{\c}$. Namely, by
conjugation this induces
$\bar\tau:\bar L_{\c}\cong \bar{\bar L_{\c}}\cong L_{\c}$, and the
composite of these two gives the identity of $L_{\c}$ (and not just an
isomorphism of
$L_{\c}$ to itself). Thus we have to choose a specific 
isomorphism $\tau:M\cong \bar M$.
If $X$ has a real  point $P$, then on the fiber over  $P$ we can choose
$\tau$ to be conjugation (and not just some constant times conjugation).

Once things are set up right, the real sections of $M$ are those
sections $s$ of $M$ such that $\tau(s)=\bar s$.
\medskip

(\ref{comp.bir.inv}) Use the fact that a birational map between
projective varieties is defined outside a codimension 2 subset.
\medskip

(\ref{gen.orient.ex}) Let $S^1\sim A\subset X(\r)$ be any loop.
Perturb it to achieve that $A$ intersects $D(\r)$ transversally at
smooth points and is disjoit from $D'(\r)$. 
\medskip

(\ref{cb.isom.ex}) $z$ is transformed by the inverse of the matrix.
Then do the explicit computation.
\medskip

(\ref{geom.cb.ex}) Let $F\to \p^1$ be a minimal conic bundle. 
There are 2 types of elementary transformations: blow up a real point
in a fiber and then contract the fiber, or  blow up  conjugate points
in conjugate fibers and then contract the fibers.

Pick any
section $C$ over $\c$.  Using elementary transformations get to the
situation when $C$ and $\bar C$ are disjoint. Show that $(C^2)=(\bar
C^2)=-m$ if there are $2m$ singular fibers. The normal form is an
affine piece of representing $F$ as a conic bundle in $\proj
f_*\o_F(C+\bar C)$. 
\medskip

(\ref{(2,2)-curves})
 The case when $D(\r)=\emptyset$ is easy. 
There are many ways to study the remaining cases.
For instance, pick a point  $P\in D(\r)$ and blow it up. By contracting
the birational transforms of the two sections through $P$, we obtain a
correspondance between pairs $(P\in D(\r))$ and triplets $(Q_1,Q_2\in
E(\r))$ where $E\subset \p^2$ is an elliptic curve. 
One has to be a little  careful since $\r\p^2$ is not orientable.
If we fix the orientations of the two copies of $\r\p^1$
in $\r\p^1\times \r\p^1$ then they  give
local orientations of $\r\p^2$ at the points $Q_1,Q_2$. There is an
ambiguity of changing both orientations (since this does not change the
orientation of $\r\p^1\times \r\p^1$).

We have to study various cases according to the location of $Q_1,Q_2$ on
$E(\r)$.  Moreover, we have to see   how the local orientations match up
if we move from $Q_1$ to $Q_2$ along $E(\r)$. 
If $Q_1,Q_2$ are both on a pseudo line, then the two local orientations
are consistent if we move in one direction and  
inconsistent in the other direction. However, if 
$Q_1,Q_2$ are both on an oval, then the two local orientations
are either  consistent in both directions or  
inconsistent if both directions. This gives 2 cases. 
\medskip

(\ref{deg4.cb.ex}) Given a  ramification curve  $D\subset \p^1\times
\p^1$  we get 2 different degree 4 Del Pezzo surfaces.  We want one
surface $S^+$ where none of the lines are real. Then, in the other
surface $S^-$, all  lines are real.  One can construct $S^-$ as
follows. Start with $\p^1\times
\p^1$, points $P_1,\dots,P_4$ in the first factor and 
$P'_1,\dots,P'_4$ in the second factor. Blow up the 4 points
$(P_i,P'_i)$.  Show that we get a degree 4 Del Pezzo surface iff there
is no isomorphism $\tau:\p^1\times \p^1$ such that $\tau(P_i)=P'_i$. 

Projecting to the first factor gives a conic fibration with singular
fibers over $P_i$. Projection to the second factor is not the right
thing to do. Instead, the second conic fibration is given by the linear
system of curves of type $(2,2)$ passing through the 4 points
$(P_i,P'_i)$.
Also keep in mind that we have to take care not only of the 4 points
but also the set $I(f)$.
\medskip 

(\ref{cubic.lines.ex}) and (\ref{deg.4.bitang}) can both be seen from
the classification.    One has to prove that we can not blow up a point
on a line. \cite{Schlafli1863} proved first that a cubic can be written
as $C_1-C_2=0$ where $C_i$ are products of linear factors and then
studied the various cases when the linear factors are all real or some
are conjugates. \cite{Zeuthen1874} notes that 2 ovals have 4 tangents,
thus we need to show that there are 4 more which are either tangents to
the same oval or at complex points. He does this by a continuity
argument. This is a bit tricky since these  4 tangents are not
invariant under deformations of the curve, just their number is.

\vskip1cm

University of Utah, Salt Lake City UT 84112 

\begin{verbatim}kollar@math.utah.edu\end{verbatim}


\medskip





 

\begin{thebibliography}{Chislenko88}

\bibitem[BCR87]{BCR87}   J. Bochnak, M. Coste and M-F. Roy,
G\'eom\'etrie alg\'ebrique r\'eelle, Springer 1987

\bibitem[Comessatti13]{Comessatti13}  A. Comessatti,
Fondamenti per la geometria sopra superfizie razionali dal punto di
vista reale, Math.   Ann.  73 (1913)  1-72

\bibitem[Comessatti14]{Comessatti14}  A. Comessatti,
Sulla connessione delle superfizie razionali  reali, Annali di Math.
23(3) (1914)  215-283

\bibitem[Hartshorne77]{Hartshorne77}   R.   Hartshorne,
Algebraic Geometry,  Springer, 1977

\bibitem[Iskovskikh67]{Iskovskikh67}    V. A. Iskovskikh,
Rational surfaces with a pencil of rational curves,
  Math. USSR Sb.   3 (1967) 563-587

\bibitem[Iskovskikh80]{Iskovskikh80}  V. A. Iskovskikh,
Minimal models of rational surfaces over arbitrary fields,  Math.
USSR Izv. 14  (1980) 17-39

\bibitem[Klein1876]{Klein1876} F. Klein, \"Uber eine neuer Art von
Riemannschen Fl\"achen, Math. Ann. 10 (1876) 398--416,  Reprinted in :
F.  Klein, Gesammelte Mathematische Abhandlungen, Springer, 1922, vol.
II.

\bibitem[Koll\'ar96]{koll96}     J.   Koll\'ar,  
  Rational Curves on Algebraic Varieties, Springer Verlag,
Ergebnisse der Math. vol. 32,  1996 

\bibitem[Koll\'ar97]{rat}  J.   Koll\'ar,
Real Algebraic Threefolds I--IV.  (in preparation)

\bibitem[Koll\'ar-Mori98]{km98}  J.   Koll\'ar - S.   Mori,
Birational geometry of algebraic varieties, Cambridge Univ. Press, 1998
(to appear)

\bibitem[Koll\'ar-Smith97]{KS97}    J.   Koll\'ar  and K. Smith, 
Rational and Non-rational Algebraic Varieties (e-prints:
alg-geom/9707013)

\bibitem[Manin66]{Manin66}    Yu. I. Manin, Rational surfaces over
perfect fields (in Russian), 
Publ. Math. IHES   30 (1966) 55-114

\bibitem[Manin72]{Manin72} Yu. I. Manin, Cubic forms (in Russian), 
 Nauka,  1972.  English translation:  North-Holland, 1974,  
 second expanded edition, 1986 
  


\bibitem[Milnor-Stasheff74]{milnor-s74} J. Milnor  and J. Stasheff,
Characteristic Classes, Princeton Univ. Press, 1974

\bibitem[Mori82]{Mori82}   S.   Mori, Threefolds whose
Canonical Bundles are not Numerically Effective, Ann.   of Math.  
 116 (1982)  133-176

\bibitem[Nash52]{Nash52} J. Nash,  Real algebraic manifolds, Ann.
Math. 56 (1952) 405-421

\bibitem[Pl\"ucker1839]{Plucker1839} J. Pl\"ucker, Theorie der
algebraischen Kurven, Bonn, 1839

\bibitem[Schl\"afli1863]{Schlafli1863} L. Schl\"afli,  On the
distribution of surfaces of the third order into species,  Phil. Trans.
Roy. Soc. London, 153(1863)193-241. Reprinted in : L. Schl\"afli,
Gesammelte Mathematische Abhandlungen, Birkh\"auser, 1953, vol. II.

\bibitem[Segre42]{Segre42} B. Segre, The non-singular cubic surfaces,
Clarendon Press, 1942

\bibitem[Segre51]{Segre51} 	B. Segre, The rational solutions of 
homogeneous cubic equations in four variables,  Notae Univ. Rosario 
  2 (1951)  1-68

\bibitem[Shafarevich72]{Shafarevich72} R. I. Shafarevich, Basic
Algebraic Geometry (in Russian), Nauka, 1972.  English translation:
Springer,  1977,   second expanded edition, 1994 

\bibitem[Silhol84]{Silhol84} R. Silhol,  Real algebraic surfaces with
rational or elliptic fibering, Math. Zeitschr. 186(1984) 465-499

\bibitem[Silhol89]{Silhol89} R. Silhol,  Real algebraic surfaces,
Springer Lecture Notes vol. 1392, 1989

\bibitem[Viro90]{Viro90} O. Ya. Viro,   Real algebraic plane curves,
Leningrad Math. J. 1 (1990) 1059-1134

\bibitem[Zeuthen1874]{Zeuthen1874} H.G. Zeuthen, Sur les diff\'erentes
formes des courbes du quatri\`eme ordre, Math. Ann. 7(1874) 410-432


\end{thebibliography}
\end{document}